\documentstyle[leqno]{article}

\newtheorem{th}{Theorem}[section]
\newtheorem{prop}[th]{Proposition}

\newcounter{defin}[section]
\renewcommand{\thedefin}{\thesection.\arabic{defin}}
\newcounter{ex}[section]

\newcounter{rem}[section]
\renewcommand{\therem}{\thesection.\arabic{rem}}
\title{Generalized Metrical Multi-Time Lagrange Model \\
for General Relativity and Electromagnetism}
\date{}
\author{Mircea Neagu}
\begin{document}
\maketitle
\begin{abstract}
Section 1 contains some physical and geometrical aspects that motivates us to
study the generalized  metrical multi-time Lagrange space of General Relativity,
denoted by $GRGML^n_p$,  whose
vertical fundamental metrical d-tensor is
$$G^{(\alpha)(\beta)}_{(i)(j)}(t^\gamma,x^k,x^k_\gamma)=h^{\alpha\beta}
(t^\gamma)e^{2\sigma(t^\gamma,x^k,x^k_\gamma)}\varphi_{ij}(x^k).
$$
Section 2  developes the geometry  of this space, in the
sense of d-connections, d-torsions and d-curvatures. Section 3  constructs the
Einstein equations of gravitational potentials of this  generalized metrical
multi-time Lagrange space. The conservation laws of the stress-energy d-tensor
of $GRGML^n_p$  are also described. Section 4 describes the Maxwell equations which govern the
electromagnetic field of this space.
\end{abstract}

\section{Geometrical and physical aspects}

\hspace{5mm} In this century, a lot of geometrical models for gravitational
and electromagnetic theories was created. We reffer especially to the well
known  Riemannian, Finslerian or, more general, Lagrangian theories.

The usual  point of
view that  the underling  geometry of space-time is  Riemannian is thought of
to be further strengthened by the constructive-axiomatic formulation
({\it EPS  conditions}) of General Relativity due to Ehlers, Pirani and Schild
\cite{4}. Within this scheme the geometry of space-time is regarded in terms of
its main substructures: conformal  and projective structures, which are in turn
thought  of to be fixed by light propagation and  freely falling non-rotating
neutral test particules, respectively.

In Finslerian context, R. Tavakol and Van  den Berg \cite{21} have showed that
the geometrical framework within which teories of gravity are sought can be
generalized without at the same  time contradicting the EPS conditions.

More general, in Lagrangian terms, a natural geometrical-axiomatic approach
of the EPS conditions is given by Miron and Anastasiei \cite{7}. In order to
discuss the EPS conditions, in Lagrangian terminology, we recall that a
{\it generalized Lagrange space} $GL^n=(M,g_{ij}(x,y))$ is defined as a pair
which consists of a real, $n$-dimensional manifold $M$ coordinated by $x=(x^i
)_{i=\overline{1,n}}$ and a {\it fundamental metrical d-tensor  $g_{ij}(x,y)$}
on $TM$, of rank $n$ and having a constante signature on $TM\backslash\{0\}$.
We point  out that $g_{ij}(x,y)$ is not necessarily 0-homogenous with respect
to the direction $y=(y^i)_{i=\overline{1,n}}$.

Let us assume that the generalized Lagrange space $GL^n=(M^n,g_{ij}(x,y))$
satisfies the following axioms:

{\bf a.1} The fundamental tensor field $g_{ij}(x,y)$ is of the form
\begin{equation}
g_{ij}(x,y)=e^{2\sigma(x,y)}\varphi_{ij}(x).
\end{equation}

{\bf a.2} The space $GL^n$ is endowed with the non-linear connection
\begin{equation}
N^i_j(x,y)=\gamma^i_{jk}(x)y^k,
\end{equation}
where $\gamma^i_{jk}(x)$ are the Christoffel symbols for the semi-Riemannian
metric $\varphi_{ij}(x)$.

The axiom {\bf a.1} asserts that the metrical d-tensor $g_{ij}(x,y)$ is conformally
to the semi-Riemannian metric $\varphi_{ij}(x)$. Therefore the spaces $GL^n$
and $R^n=(M,\varphi_{ij})$ have  the same  causal structure.

The axiom {\bf a.2} enssures us that the  autoparallel curves of the nonlinear
connection $N^i_j(x,y)$ of the space $GL^n$ coincide to  the geodesics of
$R^n$.

Consequently, under above axiomatic assumptions, the generalized Lagrange
space $GL^n$ becomes  a convenient mathematical  model for General  Relativity,
because it verifies the EPS conditions. The differential geometry of the
generalized Lagrange  space $GL^n=(M,e^{2\sigma(x,y)}\varphi_{ij}(x))$ is now
completely developed in \cite{1}, \cite{7}, \cite{9}.

It is well known that the jet fibre
bundle of order one $J^1(T,M)$ is a basic object in the study of classical and quantum
field theories.
In a general setting, in  a previous  paper  \cite{10}, Neagu creates a natural
geometry  of physical fields induced by a Kronecker $h$-regular vertical
metrical d-tensor $G^{(\alpha)(\beta)}_{(i)(j)}(t^\gamma,x^k,x^k_\gamma)$
on the total space  of the 1-jet
vector bundle $J^1(T,M)\to T\times M$, where $(T, h)$  is a smooth, real,
$p$-dimensional semi-Riemannian manifold coordinated by $t=(t^\alpha)_{\alpha
=\overline{1,p}}$, whose physical meaning is that of {\it "multidimensional
time"}. We point out that $J^1(T,M)$ is coordinated by $(t^\alpha,x^i,x^i_
\alpha)$, where $x^i_\alpha$ have the physical meaning of {\it partial
directions}.

In the multi-temporal context, the fundamental geometrical
concept used in the geometrization of a vertical multi-time metrical d-tensor
$G^{(\alpha)(\beta)}_{(i)(j)}(t^\gamma,x^k,x^k_\gamma)$
is that of {\it generalized metrical multi-time Lagrange space} \cite{10}.
This geometrical
concept with physical meaning is
represented by a pair $GML^n_p=(J^1(T,M),G^{(\alpha)(\beta)}_{(i)(j)})$
consisting of the 1-jet space and a {\it Kronecker $h$-regular} vertical
multi-time metrical d-tensor $G^{(\alpha)\beta)}_{(i)(j)}$, that  is
\begin{equation}
G^{(\alpha)(\beta)}_{(i)(j)}(t^\gamma,x^k,x^k_\gamma)=h^{\alpha\beta}(t^
\gamma)g_{ij}(t^\gamma,x^k,x^k_\gamma),
\end{equation}
where $g_{ij}(t^\gamma,x^k,x^k_\gamma)$ is a d-tensor on $J^1(T,M)$,
symmetric, of rank $n$ and having a constant signature. The d-tensor
$g_{ij}(t^\gamma,x^k,x^k_\gamma)$ is called the
{\it spatial metrical d-tensor of $GML^n_p$}.

The differential geometry of the generalized metrical multi-time Lagrange
spaces together with its attached field theory are now
considerably developed in \cite{10}.

Following the general physical and geometrical development from \cite{10},
the aim of this paper is to study the generalized metrical multi-time
Lagrange space $GRGML^n_p$, whose spatial metrical d-tensor is of the form
\begin{equation}
g_{ij}(t^\gamma,x^k,x^k_\gamma)=e^{2\sigma(t^\gamma,x^k,x^k_\gamma)}\varphi_
{ij}(x^k),
\end{equation}
where $\varphi_{ij}(x^k)$  is  a semi-Riemannian metric on the spatial
manifold  $M$ and \linebreak $\sigma:J^1(T,M)\to R$ is a conformal smooth
function,  which gives the {\it magnitude} of directions $x^i_\alpha$.
\medskip\\
\addtocounter{rem}{1}
{\bf Remark \therem} From  physical point of view,  the  interesting properties
of this space are obtained  considering the special conformal  functions:
\medskip

i) $\sigma=U^{(\alpha)}_{(i)}(t^\gamma,x^k)x^i_\alpha$,\medskip

ii) $\sigma=h^{\alpha\beta}(t^\gamma)A_i(x^k)A_j(x^k)x^i_\alpha x^j_\beta$,
\medskip

iii) $\sigma=\varphi_{ij}(x^k)X^\alpha(t^\gamma)X^\beta(t^\gamma)x^i_\alpha
x^j_\beta$,\medskip\\
where $U^{(\alpha)}_{(i)}(t^\gamma,x^k)$ is a d-tensor on $E$, $A_i(x^k)$ is
a covector field on $M$, and $X^\alpha(t^\gamma)$ is a vector field on $T$.
For more details, see \cite{7}.
\medskip

In order to develope the geometry of this generalized metrical multi-time
Lagrange space, we need  a nonlinear connection $\Gamma=(M^{(i)}_{(\alpha)\beta},
N^{(i)}_{(\alpha)j})$ on $J^1(T,M)$. In this direction, we fix {\it "a priori"}
the nonlinear  connection $\Gamma$ defined by the temporal components
\begin{equation}
M^{(i)}_{(\alpha)\beta}=-H^\mu_{\alpha\beta}x^i_\mu
\end{equation}
and the spatial components
\begin{equation}
N^{(i)}_{(\alpha)j}=\gamma^i_{jm}x^m_\alpha,
\end{equation}
where $H^\alpha_{\beta\gamma}$ (resp. $\gamma^i_{jk}$) are the Christoffel
symbols  of the semi-Riemannian metric $h_{\alpha\beta}$ (resp. $\varphi_{ij}$).
\medskip\\
\addtocounter{rem}{1}
{\bf Remarks \therem} i) The previous nonlinear connection $\Gamma$  is dependent
only the vertical fundamental metrical d-tensor $G^{(\alpha)(\beta)}_{(i)(j)}$
of $GRGML^n_p$. This fact emphasize the {\it metrical character} of the geometry
attached to this space, i. e. , all geometrical  objects are directly arised from
$G^{(\alpha)(\beta)}_{(i)(j)}$.

ii) The spatial components $N^{(i)}_{(\alpha)j}$ of the fixed nonlinear
connection $\Gamma$ are  {\it without torsion} \cite{10}.\medskip

Investigating the possibility of realising the {\it multi-time EPS conditions}
for General Relativity, in a general setting, let us start with the generalized
metrical multi-time Lagrange space
\begin{equation}
GRGML^n_p=(J^1(T,M),G^{(\alpha)(\beta)}_{(i)(j)}=h^{\alpha\beta}(t^\gamma)
g_{ij}(t^\gamma,x^k,x^k_\gamma)),
\end{equation}
which verifies the axioms {\bf A.1} and {\bf A.2} from below:
\medskip

{\bf A.1} The vertical  fundamental tensor field $G^{(\alpha)(\beta)}_{(i)(j)}
(t^\gamma,x^k,x^k_\gamma)$ is of the form
\begin{equation}\label{A.1}
G^{(\alpha)(\beta)}_{(i)(j)}(t^\gamma,x^k,
x^k_\gamma)=h^{\alpha\beta}(t^\gamma)e^{2\sigma(t^\gamma,x^k,x^k_\gamma)}
\varphi_{ij}(x^k).
\end{equation}

{\bf A.2} The generalized  metrical multi-time Lagrange space $GRGML^n_p$ is
endowed with the non-linear connection
$\Gamma=(M^{(i)}_{(\alpha)\beta},N^{(i)}_{(\alpha)j})$, defined  by the
components
\begin{equation}\label{A.2}
M^{(i)}_{(\alpha)\beta}=-H^\mu_{\alpha\beta}x^i_\mu,
\qquad
N^{(i)}_{(\alpha)j}=\gamma^i_{jm}x^m_\alpha,
\end{equation}
where $H^\alpha_{\beta\gamma}$ (resp. $\gamma^i_{jk}$) are the Christoffel
symbols  of the semi-Riemannian metric $h_{\alpha\beta}$ (resp. $\varphi_{ij}$).
\medskip\\

Let us consider the Lagrangian function $L:J^1(T,M)\to R$, used in the Polyakov
model of bosonic strings,
\begin{equation}
L(t^\gamma,x^k,x^k_\gamma)=h^{\alpha\beta}(t^\gamma)\varphi_{ij}(x^k)
x^i_\alpha x^j_\beta,
\end{equation}
and its vertical fundamental metrical d-tensor,
\begin{equation}
g^{(\alpha)(\beta)}_{(i)(j)}={1\over 2}{\partial^2L\over\partial x^i_\alpha
\partial x^j_\beta}=h^{\alpha\beta}(t^\gamma)\varphi_{ij}(x^k).
\end{equation}
\addtocounter{rem}{1}
{\bf Remark \therem} The extremals of the Lagrangian ${\cal L}=L\sqrt{\vert
h\vert}$ are exactly the harmonic maps between the semi-Riemannian spaces
$(T,h)$ and $(M,\varphi)$ \cite{3}.\medskip

Let $EDML^n_p=(J^1(T,M),L)$ be  the {\it metrical multi-time Lagrange space of
electrodynamics} corresponding to $L$ \cite{11}, \cite{15}.
In this context, we have the following important result
\begin{th}
i) The generalized metrical multi-time Lagrange space $GRGML^n_p$ has the same
conformal structure as the autonomous metrical multi-time  Lagrange space
of electrodynamics $EDML^n_p=(J^1(T,M),L)$.

ii) The harmonic maps of the nonlinear connection $\Gamma$ of $GRGML^n_p$
coincide with those of the canonical nonlinear connection \cite{14} of the
autonomous metrical multi-time Lagrange  space of electrodynamics
$EDML^n_p=(J^1(T,M),L)$. Moreover, these are exactly the harmonic maps between
the semi-Riemannian spaces $(T,h)$ and $(M,\varphi)$.
\end{th}
{\bf Proof.}  The axiom {\bf A.1} enssures us that both spaces $GRGML^n_p$ and
$EDML_p^n$ have the same causal  structure.

Following  the paper \cite{11}, we deduce that  the components of the  canonical
nonlinear connection of $EDML^n_p$ are described exactly by formulas \ref{A.2},
that is, those from
the axiom {\bf  A.2}. Moreover, the relation between sprays and the
components of a nonlinear connection, described in \cite{14}, and the
definition of harmonic maps attached to a given multi-time dependent spray on
$J^1(T,M)$, imply what we were looking for. \rule{5pt}{5pt}\medskip

In conclusion, we can assert that the generalized metrical multi-time Lagrange
space $GRGML^n_p$, which verifies the assumptions {\bf A.1} and {\bf A.2},
represents a convenient relativistic model, in the  multi-temporal context,
since it has the same conformal and projective properties as the  autonomous
metrical multi-time Lagrange space from Polyakov model of bosonic strings.

\section{Cartan canonical connection}

\setcounter{equation}{0}
\hspace{5mm} In this section, we will apply the general geometrical development
of the generalized metrical multi-time Lagrange spaces \cite{10}, to the
particular space
\begin{equation}\label{mdt}
GRGML^n_p=(J^1(T,M),G^{(\alpha)(\beta)}_{(i)(j)}=h^{\alpha\beta}(t^\gamma)e^{2
\sigma(t^\gamma,x^k,x^k_\gamma)}\varphi_{ij}(x^k)),
\end{equation}
endowed with the nonlinear  connection $\Gamma=(M^{(i)}_{(\alpha)\beta},
N^{(i)}_{(\alpha)j})$, where
\begin{equation}\label{nc}
M^{(i)}_{(\alpha)\beta}=-H^\mu_{\alpha\beta}x^i_\mu,
\qquad
N^{(i)}_{(\alpha)j}=\gamma^i_{jm}x^m_\alpha.
\end{equation}

Let $\displaystyle{\left\{{\delta\over\delta
t^\alpha}, {\delta\over\delta x^i}, {\partial\over\partial x^i_\alpha}\right\}
\subset{\cal X}(J^1(T,M))}$ and $\{dt^\alpha, dx^i, \delta x^i_\alpha\}\subset{\cal
X}^*(J^1(T,M))$  be the adapted bases of the nonlinear connection $\Gamma$, where
\begin{equation}
\left\{\begin{array}{l}\medskip
\displaystyle{{\delta\over\delta t^\alpha}={\partial\over\partial t^\alpha}-
M^{(j)}_{(\beta)\alpha}{\partial\over\partial x^j_\beta}}\\\medskip
\displaystyle{{\delta\over\delta x^i}={\partial\over\partial x^i}-
N^{(j)}_{(\beta)i}{\partial\over\partial x^j_\beta}}\\
\delta x^i_\alpha=dx^i_\alpha+M^{(i)}_{(\alpha)\beta}dt^\beta+N^{(i)}_{(\alpha)
j}dx^j.
\end{array}\right.
\end{equation}

Following the paper \cite{10}, by a direct calculation, we can  determine the {\it
Cartan canonical connection} of $GRGML^n_p$, together with its torsion and
curvature local d-tensors.
\begin{th}
The Cartan canonical connection
$
C\Gamma=(H^\gamma_{\alpha\beta}, G^k_{j\gamma}, L^i_{jk}, C^{i(\gamma)}_
{j(k)})
$
of $GRGML^n_p$ has the adapted coefficients
\begin{equation}
\begin{array}{l}\medskip
H^\gamma_{\alpha\beta}=H^\gamma_{\alpha\beta},\quad
G^k_{j\gamma}=\sigma_\gamma\delta^k_j, \quad
L^k_{ij}=\gamma^k_{ij}+\Lambda^k_{ij},\\
C^{i(\gamma)}_{j(k)}=\sigma^{(\gamma)}_{(k)}\delta^i_j+\sigma^{(\gamma)}_{(j)}
\delta^i_k-\varphi_{jk}\sigma^{i\gamma},
\end{array}
\end{equation}
where
\begin{equation}
\begin{array}{l}\medskip
\displaystyle{\sigma_\gamma={\delta\sigma\over\delta t^\gamma},\quad
\sigma_k={\delta\sigma\over\delta x^k},\quad
\sigma^{(\gamma)}_{(k)}={\partial\sigma\over\partial x^k_\gamma},}\\
\sigma^{k\gamma}=\varphi^{km}\sigma^{(\gamma)}_{(m)},\quad
\Lambda^k_{ij}=\sigma_i\delta^k_j+\sigma_j\delta^k_i-\varphi_{ij}\sigma^k.
\end{array}
\end{equation}
\end{th}

\begin{th}
The torsion {\em\bf T} of the Cartan canonical connection of $GRGML^n_p$
is determined by seven effective local d-tensors, namely,
\begin{equation}
\begin{array}{l}\medskip
T^m_{\alpha j}=-\sigma_\alpha\delta^m_j,\quad P^{m(\beta)}_{i(j)}=C^{m(\beta)}
_{i(j)},\quad P^{(m)\;\;(\beta)}_{(\mu)\alpha(j)}=-\sigma_\alpha\delta^\beta_\mu
\delta^m_j,\\\medskip
P^{(m)\;(\beta)}_{(\mu)i(j)}=-\Lambda^m_{ij}\delta^\beta_\mu,\quad
S^{(m)(\alpha)(\beta)}_{(\mu)(i)(j)}=\delta^\alpha_\mu C^{m(\beta)}_{i(j)}-
\delta^\beta_\mu C^{m(\alpha)}_{i(j)},\\
R^{(m)}_{(\mu)\alpha\beta}=-H^\gamma_{\mu\alpha\beta}x^m_\gamma,
\quad R^{(m)}_{(\mu)\alpha j}=0,\quad
R^{(m)}_{(\mu)ij}=r^m_{ijk}x^k_\mu,
\end{array}
\end{equation}
where $H^\gamma_{\mu\alpha\beta}$ (resp. $r^m_{ijk}$) are the local curvature
tensors of the semi-Riemannian metric $h_{\alpha\beta}$ (resp. $\varphi_{ij}$).
\end{th}

In order to describe  the local curvature  d-tensors of Cartan canonical
connection of $GRGML^n_p$, let us consider $B\Gamma=(H^\gamma_{\alpha\beta},0,
\gamma^i_{jk},0)$, the Berwald $h$-normal $\Gamma$-linear connection attached
to the semi-Riemannian metrics $h_{\alpha\beta}$ and $\varphi_{ij}$ \cite{13}.
We denote by $"_{//\alpha}"$, $"_{\Vert_i}"$ and $"\Vert^{(\alpha)}_{(i)}"$,
the local covariant derivatives  induced by $B\Gamma$. Now, taking into account
the expressions of these local covariant derivatives \cite{13}, by a direct
calculation, we deduce
\begin{prop}
The Berwald connection $B\Gamma$ of $GRGML^n_p$ has the following metrical
properties:
\begin{equation}
\left\{\begin{array}{lll}\medskip
h_{\alpha\beta//\gamma}=0,\quad h_{\alpha\beta\Vert k}=0,\quad
h_{\alpha\beta}\Vert^{(\gamma)}_{(k)}=0,\\\medskip
\varphi_{ij//\gamma}=0,\quad \varphi_{ij\Vert k}=0,\quad
\varphi_{ij}\Vert^{(\gamma)}_{(k)}=0,\\
g_{ij//\gamma}=2\sigma_\gamma g_{ij},\quad g_{ij\Vert k}=2\sigma_k g_{ij},
\quad  g_{ij}\Vert^{(\gamma)}_{(k)}=2\sigma^{(\gamma)}_{(k)}g_{ij},
\end{array}\right.
\end{equation}
where $g_{ij}(t^\gamma,x^k,x^k_\gamma)=e^{2\sigma(t^\gamma,x^k,x^k_\gamma)}
\varphi_{ij}(x^k)$.
\end{prop}

In these conditions, using the general expressions of  the local curvature
d-tensors attached to the Cartan canonical connection of a generalized  metrical
multi-time Lagrange space $GML^n_p$, by computations, we obtain

\begin{th}
The curvature {\em\bf R} of the Cartan canonical connection of $GRGML^n_p$
is determined by seven effective local d-tensors, expressed by,
\begin{equation}\label{curv}
\begin{array}{l}\medskip
\displaystyle{H^\alpha_{\eta\beta\gamma}={\partial H^\alpha_{\eta\beta}\over
\partial t^\gamma}-{\partial H^\alpha_{\eta\gamma}\over\partial t^\beta}+
H^\mu_{\eta\beta}H^\alpha_{\mu\gamma}-H^\mu_{\eta\gamma}H^\alpha_{\mu\beta},}
\\\medskip
\displaystyle{R^l_{i\beta\gamma}=\left[C^{l(\mu)}_{i(m)}-\delta^l_i\sigma^{(\mu)}
_{(m)}\right]R^{(m)}_{(\mu)\beta\gamma},}
\\\medskip
\displaystyle{R^l_{i\beta k}=\left[\varphi_{ik}\varphi^{lm}-\delta^l_k
\delta^m_i\right]\sigma_{m//\beta},}
\\\medskip
\displaystyle{R^l_{ijk}=r^l_{ijk}+\Lambda^l_{ij\Vert k}-\Lambda^l_{ik\Vert j}
+\Lambda^m_{ij}\Lambda^l_{mk}-\Lambda^m_{ik}\Lambda^l_{mj}+C^{l(\mu)}_{i(m)}
R^{(m)}_{(\mu)jk},}
\\\medskip
\displaystyle{P^{l\;\;(\gamma)}_{i\beta(k)}={\partial \sigma_{\beta}\over
\partial x^k_\gamma}\delta^l_i-{\delta C^{l(\gamma)}_{i(k)}\over\delta
t^\beta}-C^{l(\mu)}_{i(k)}H^\gamma_{\mu\beta},}
\\\medskip
\displaystyle{P^{l\;(\gamma)}_{ij(k)}=\Lambda^l_{ij}\Vert_{(k)}^{(\gamma)}-
C^{l(\gamma)}_{i(k)\Vert j}-\Lambda^l_{jm}C^{m(\gamma)}_{i(k)}+\Lambda^m_{ij}
C^{l(\gamma)}_{m(k)},}
\\\medskip
\displaystyle{S^{l(\beta)(\gamma)}_{i(j)(k)}=C^{l(\beta)}_{i(j)}\Vert_{(k)}^
{(\gamma)}-C^{l(\gamma)}_{i(k)}\Vert_{(j)}^{(\beta)}+C^{m(\beta)}_{i(j)}
C^{l(\gamma)}_{m(k)}-C^{m(\gamma)}_{i(k)}C^{l(\beta)}_{m(j)},}
\end{array}
\end{equation}
where  $H^\eta_{\alpha\beta\gamma}$ and $r^l_{ijk}$ are the curvature tensors
of the semi-Riemannian metrics $h_{\alpha\beta}$ and $\varphi_{ij}$.
\end{th}

In order to give a more  natural and beautiful form to the local curvature
d-tensors \ref{curv}, we need the following notations:
\begin{equation}
\begin{array}{l}\medskip
\varphi^{ij}_{kl}=\varphi^{ij}\varphi_{kl}-\delta^i_l\delta^j_k,\\\medskip
\displaystyle{\sigma_{jk}=\sigma_{j\Vert k}-\sigma_j\sigma_k+{1\over 2}
\varphi_{jk}<\sigma,\sigma>,}\\\medskip
\displaystyle{\sigma_{j(k)}^{\;\;(\gamma)}=\sigma_j\Vert^{(\gamma)}_
{(k)}-\sigma_k\sigma_{(j)}^{(\gamma)}+{1\over 2}\varphi_{jk}<\sigma,\sigma>^
\gamma,}\\\medskip
\displaystyle{\sigma_{(j)k}^{(\beta)}=\sigma_{(j)\Vert k}^
{(\beta)}-\sigma_j\sigma_{(k)}^{(\beta)}+{1\over 2}\varphi_{jk}<\sigma,\sigma>^
\beta,}\\\medskip
\displaystyle{\sigma_{(j)(k)}^{(\beta)(\gamma)}=\sigma_{(j)}^{(\beta)}
\Vert_{(k)}^{(\gamma)}-\sigma_{(j)}^{(\beta)}\sigma_{(k)}^{(\gamma)}+
{1\over 2}\varphi_{jk}<\sigma,\sigma>^{\beta\gamma},}
\end{array}
\end{equation}
where
\begin{equation}
\begin{array}{l}\medskip
<\sigma,\sigma>=\varphi^{rm}\sigma_r\sigma_m,\\\medskip
<\sigma,\sigma>^\beta=\varphi^{rm}\sigma_r\sigma_{(m)}^{(\beta)}=
\varphi^{rm}\sigma_{(r)}^{(\beta)}\sigma_m,\\\medskip
<\sigma,\sigma>^{\beta\gamma}=\varphi^{rm}\sigma_{(r)}^{(\beta)}
\sigma_{(m)}^{(\gamma)}=\varphi^{rm}\sigma_{(r)}^{(\gamma)}\sigma_{(m)}^{(
\beta)}=<\sigma,\sigma>^{\gamma\beta}.
\end{array}
\end{equation}

The Ricci identities of the Berwald connection \cite{16}, applied  to the d-tensors
$\sigma_i$ and $\sigma^{(\alpha)}_{(i)}$, imply
\begin{prop}
The following tensorial identities are true:
\begin{equation}
\begin{array}{l}\medskip
\sigma_{jk}-\sigma_{kj}=-r^m_{ljk}\sigma^{(\mu)}_{(m)}x^l_\mu,\\\medskip
\sigma_{i(j)}^{\;\;(\beta)}-\sigma_{(j)i}^{(\beta)}=0\\
\sigma_{(i)(j)}^{(\alpha)(\beta)}-\sigma_{(j)(i)}^{(\beta)(\alpha)}=0.
\end{array}
\end{equation}
\end{prop}

Consequently, by simple computations, we obtain

\begin{th}
The local curvature d-tensors \ref{curv} of the Cartan canonical connection
of $GRGML^n_p$ have  the  following new expressions:
\begin{equation}\label{ncurv}
\hspace*{5mm}
\begin{array}{l}\medskip
H^\alpha_{\eta\beta\gamma}=H^\alpha_{\eta\beta\gamma},
\\\medskip
R^l_{i\beta\gamma}=-\varphi^{lm}_{ir}H^\varepsilon_{\mu\beta\gamma}\sigma^{(
\mu)}_{(m)}x^r_\varepsilon,
\\\medskip
R^l_{i\beta k}=\varphi_{ik}^{lm}\sigma_{m//\beta},
\\\medskip
R^l_{ijk}=r^l_{ijk}+\sigma_{ik}\delta^l_j-\sigma_{ij}\delta^l_k+\varphi^{ls}
\left[\varphi_{ik}\sigma_{sj}-\varphi_{ij}\sigma_{sk}\right]-\varphi^{ml}_{si}
r^s_{pjk}\sigma^{(\mu)}_{(m)}x^p_\mu,
\\\medskip
P^{l\;\;(\gamma)}_{i\beta(k)}=\varphi^{lm}_{ik}\sigma_{(m)//\beta}^
{(\gamma)}+\sigma_\beta\sigma^{(\gamma)}_{(k)}\delta^l_i,
\\\medskip
P^{l\;(\gamma)}_{ij(k)}=\sigma_{i(k)}^{\;\;(\gamma)}\delta^l_j-\sigma_{(i)j}^
{(\gamma)}\delta^l_k-\varphi^{ls}\left[\varphi_{ij}\sigma_{s(k)}^{\;\;(\gamma)}
-\varphi_{ik}\sigma_{(s)j}^{(\gamma)}\right]+\\\medskip
\hspace*{13mm}+\varphi_{jk}\left[\sigma^l\sigma^{(\gamma)}_{(i)}-\sigma_i\sigma^{l\gamma}\right],
\\\medskip
S^{l(\beta)(\gamma)}_{i(j)(k)}=\sigma_{(i)(k)}^{(\beta)(\gamma)}\delta^l_j-
\sigma_{(i)(j)}^{(\gamma)(\beta)}\delta^l_k-\varphi^{ls}\left[\varphi_{ij}
\sigma_{(s)(k)}^{(\beta)(\gamma)}-\varphi_{ik}\sigma_{(s)(j)}^{(\gamma)
(\beta)}\right]+\\
\hspace*{16mm}+\varphi_{jk}\left[\sigma^{l\beta}\sigma^{(\gamma)}_{(i)}-
\sigma^{l\gamma}\sigma^{(\beta)}_{(i)}\right],
\end{array}
\end{equation}
where  $\sigma^l=\varphi^{lm}\sigma_m$ and $\sigma^{l\mu}=\varphi^{lm}
\sigma_{(m)}^{(\mu)}$.
\end{th}

\section{Einstein equations of gravitational field}

\setcounter{equation}{0}

\hspace{5mm} In order to develope the gravitational theory on $GRGML^n_p$, we
point out that the vertical metrical d-tensor \ref{mdt} and its fixed
nonlinear connection \ref{nc} induce a natural {\it gravitational
$h$-potential} on the 1-jet space $J^1(T,M)$ (i.  e. a Sasakian-like metric),
which is expressed by \cite{10}
\begin{equation}
G=h_{\alpha\beta}dt^\alpha\otimes dt^\beta+e^{2\sigma}\varphi_{ij}dx^i\otimes
dx^j+h^{\alpha\beta}e^{2\sigma}\varphi_{ij}\delta x^i_\alpha\otimes\delta x^j
_\beta.
\end{equation}
Let $C\Gamma=(H^\gamma_{\alpha\beta},G^k_{j\gamma},L^i_{jk},C^{i(\gamma)}_
{j(k)})$ be the Cartan canonical connection of $GRGML^n_p$.

We postulate that the Einstein equations which govern the gravitational $h$-potential
$G$ of $GRGML^n_p$ are the Einstein equations attached to the Cartan canonical
connection and the adapted metric $G$ on $J^1(T,M)$, that is,
\begin{equation}
Ric(C\Gamma)-{Sc(C\Gamma)\over 2}G={\cal K}{\cal T},
\end{equation}
where $Ric(C\Gamma)$ represents the Ricci d-tensor of the Cartan connection,
$Sc(C\Gamma)$ is its scalar curvature, ${\cal K}$ is the Einstein constant
and ${\cal T}$ is an intrinsec distinguished tensor of matter which is called  the {\it
stress-energy} d-tensor.

In the adapted basis $(X_A)=\displaystyle{\left({\delta\over\delta t^\alpha},
{\delta\over\delta x^i},{\partial\over\partial x^i_\alpha}\right)}$ attached
to $\Gamma$, the curvature d-tensor {\bf R} of the Cartan connection is
expressed locally by {\bf R}$(X_C,X_B)X_A=R^D_{ABC}X_D$. Hence, it follows
that we have $R_{AB}=Ric(C\Gamma)(X_A,X_B)=R^D_{ABD}$ and $Sc(C\Gamma)=G^{AB}
R_{AB}$, where
\begin{equation}
G^{AB}=\left\{\begin{array}{ll}\medskip
h_{\alpha\beta},&\mbox{for}\;\;A=\alpha,\;B=\beta\\\medskip
e^{-2\sigma}\varphi^{ij},&\mbox{for}\;\;A=i,\;B=j\\\medskip
h_{\alpha\beta}e^{-2\sigma}\varphi^{ij},&\mbox{for}\;\;A={(i)\atop(\alpha)},\;B={(j)\atop(\beta)}\\
0,&\mbox{otherwise}.
\end{array}\right.
\end{equation}

Taking into account the expressions \ref{ncurv} of the local curvature
d-tensors of the Cartan connection  of $GRGML^n_p$, we obtain without difficulties
\begin{th}
The Ricci  d-tensor $Ric(C\Gamma)$ of $GRGML^n_p$ is determined by seven effective
local d-tensors, expressed, in adapted basis, by:
\begin{equation}
\hspace*{5mm}
\begin{array}{l}\medskip
R_{\alpha\beta}=H_{\alpha\beta}=H^\mu_{\alpha\beta\mu}
\\\medskip
R_{i\beta}=(1-n)\sigma_{i//\beta},
\\\medskip
R_{ij}=r_{ij}+(2-n)\sigma_{ij}-\varphi_{ij}<\sigma>
+r_{mj}\sigma^{(\mu)}_{(i)}x^m_\mu-\varphi_{is}r^s_{mjp}\sigma^{p\mu}x^m_\mu,
\\\medskip
P^{(\alpha)}_{(i)\beta}=(1-n)\sigma_{(i)//\beta}^{(\alpha)}+\sigma^{(\alpha)}
_{(i)}\sigma_\beta,
\\\medskip
P^{(\alpha)}_{(i)j}=\sigma_{i(j)}^{\;(\alpha)}+(1-n)\sigma_{(i)j}^{(\alpha)}
-\varphi_{ij}<\sigma>^\alpha+\sigma_j\sigma_{(i)}^{(\alpha)}-\sigma_i\sigma^
{(\alpha)}_{(j)},
\\\medskip
P^{\;(\alpha)}_{i(j)}=P^{(\alpha)}_{(i)j},
\\\medskip
S^{(\beta)(\gamma)}_{(j)(k)}=(1-n)\sigma_{(j)(k)}^{(\beta)(\gamma)}+
\sigma_{(j)(k)}^{(\gamma)(\beta)}-\varphi_{jk}<\sigma>^{\beta\gamma}+
\sigma^{(\beta)}_{(j)}\sigma^{(\gamma)}_{(k)}-\sigma^{(\beta)}_{(k)}\sigma^
{(\gamma)}_{(j)},
\end{array}
\end{equation}
where $n=\dim M$, $H_{\alpha\beta}$ (resp.
$r_{ij}$) are the local Ricci tensors of the metric
$h_{\alpha\beta}$ (resp. $\varphi_{ij}$),
$\quad<\sigma>=\varphi^{rs}\sigma_{rs}$,$\quad<\sigma>^\alpha=\varphi^{rs}\sigma_
{r(s)}^{\;(\alpha)}$ and $<\sigma>^{\alpha\beta}=\varphi^{rs}\sigma_{(r)(s)}
^{(\alpha)(\beta)}$.
\end{th}

Let us denote $H=h^{\alpha\beta}H_{\alpha\beta}$, $R=e^{2\sigma}\varphi^{ij}
R_{ij}$ and $S=h_{\alpha\beta}e^{2\sigma}\varphi^{ij}S^{(\alpha)(\beta)}_
{(i)(j)}$. In this context, by a simple  calculation, it follows
\begin{th}
The scalar curvature of the Cartan connection $C\Gamma$ of $GRGML^n_p$ is
given by
\begin{equation}
Sc(C\Gamma)=H+R+S,
\end{equation}
where
\begin{equation}
\begin{array}{l}\medskip
H=h^{\alpha\beta}H_{\alpha\beta},\\\medskip
R=e^{-2\sigma}\left[r+2(1-n)<\sigma>+2r_{ms}\sigma^{s\mu}x^m_\mu\right],\\
S=2(1-n)e^{-2\sigma}<<\sigma>>,
\end{array}
\end{equation}
where $n=\dim M$, $H$ (resp $r$) is the scalar curvature of the semi-Riemannian
metric $h_{\alpha\beta}$ (resp. $\varphi_{ij}$) and $<<\sigma>>=h_{\alpha\beta}<\sigma>^{\alpha\beta}$.
\end{th}

Following the gravitational field theoretical exposition on a generalized metrical
multi-time Lagrange space $GML^n_p$ from the paper \cite{10}, by  local computations,
we can give
\begin{th}
If $p>2$ and $n>2$, the Einstein equations which govern the gravitational
$h$-potential $G$ of $GRGML^n_p$ are
$$
\left\{\begin{array}{l}\medskip
\displaystyle{H_{\alpha\beta}-{H\over 2}h_{\alpha\beta}={\cal K}\tilde{\cal T}_
{\alpha\beta}}\\\medskip
\displaystyle{r_{ij}-{r\over 2}\varphi_{ij}+\rho_{ij}={\cal K}\tilde{\cal T}_{ij}}\\
\displaystyle{S^{(\alpha)(\beta)}_{(i)(j)}+(n-1)<<\sigma>>h^{\alpha\beta}
\varphi_{ij}={\cal K}\tilde{\cal T}^{(\alpha)(\beta)}_{(i)(j)}},
\end{array}\right.\leqno{(E^\prime_1)}
$$
$$
\left\{\begin{array}{lll}\medskip
0={\cal T}_{\alpha i},&R_{i\alpha}={\cal K}{\cal T}_{i\alpha},&
P^{(\alpha)}_{(i)\beta}={\cal K}{\cal T}^{(\alpha)}_{(i)\beta}\\
0={\cal T}^{\;(\beta)}_{\alpha(i)},&
P^{\;(\alpha)}_{i(j)}={\cal K}{\cal T}^{\;(\alpha)}_{i(j)},&
P^{(\alpha)}_{(i)j}={\cal K}{\cal T}^{(\alpha)}_{(i)j},
\end{array}\right.\leqno{(E_2)}
$$
where $\tilde{\cal T}_{\alpha\beta}$, $\tilde{\cal T}_{ij}$ and $\tilde{\cal
T}^{(\alpha)(\beta)}_{(i)(j)}$  represent the components  of the new
stress-energy d-tensor $\tilde{\cal T}$, expressed by
\begin{equation}\label{*}
\left\{\begin{array}{l}\medskip
\displaystyle{\tilde{\cal T}_{\alpha\beta}={\cal T}_{\alpha\beta}+{R+S\over
2{\cal K}}h_{\alpha\beta}}\\\medskip
\displaystyle{\tilde{\cal T}_{ij}={\cal T}_{ij}+{H+S\over 2{\cal K}}
e^{2\sigma}\varphi_{ij}}\\
\displaystyle{\tilde{\cal T}^{(\alpha)(\beta)}_{(i)(j)}={\cal T}^{(\alpha)
(\beta)}_{(i)(j)}+{H+R\over 2{\cal K}}h^{\alpha\beta}e^{2\sigma}\varphi_{ij}},
\end{array}\right.
\end{equation}
and
\begin{equation}\hspace*{8mm}
\rho_{ij}=(2-n)\sigma_{ij}-(3-2n)<\sigma>\varphi_{ij}-[\varphi_{is}r^s_{mjp}+
2r_{mp}\varphi_{ij}-r_{mj}\varphi_{ip}]\sigma^{p\mu}x^m_\mu.
\end{equation}
\end{th}
\addtocounter{rem}{1}
{\bf Remarks \therem}  i) It is remarkable that, in the particular case $\sigma
\equiv 0$, the Einstein equations $(E_1^\prime)$ of $GRGML^n_p$ reduce to the classical ones.

ii) Note that, in order to have the compatibility of the Einstein equations, it is
necessary that the certain  adapted local components of the stress-energy
d-tensor vanish {\it "a priori"}.
\medskip

From physical point of view, the
stress-energy d-tensor ${\cal T}$ must  verify the local {\it conservation
laws} ${\cal T}^B_{A\vert B}=0,\;\forall\;A\in\{\alpha,i,{(\alpha)\atop (i)}\}$,
where ${\cal T}^B_A=G^{BD}{\cal T}_{DA}$, $"_{\vert A}"$, represent one of the
local covariant derivatives $"_{/\beta}"$, $"_{\vert j}"$ or $"\vert^{(\beta)}
_{(j)}"$, associated to the Cartan canonical connection $C\Gamma$ \cite{13}.

In this context, let us denote
\begin{equation}
\begin{array}{lll}\medskip
\tilde{\cal T}_T=h^{\alpha\beta}\tilde{\cal T}_{\alpha\beta},&
\tilde{\cal T}_M=e^{-2\sigma}\varphi^{ij}\tilde{\cal T}_{ij},&
\tilde{\cal T}_v=h_{\mu\nu}e^{-2\sigma}\varphi^{mr}\tilde{\cal
T}^{(\mu)(\nu)}_{(m)(r)},\\
\tilde{\cal T}^\alpha_\beta=h^{\alpha\mu}\tilde{\cal T}_{\mu\beta},&
\tilde{\cal T}^i_j=e^{-2\sigma}\varphi^{im}\tilde{\cal T}_{mj},&
\tilde{\cal T}^{(i)(\beta)}_{(\alpha)j}=h_{\alpha\mu}e^{-2\sigma}\varphi^{mi}
\tilde{\cal T}^{(\mu)(\beta)}_{(m)(i)}.
\end{array}
\end{equation}

Following again the development of gravitational generalized metrical
multi-time theory from \cite{10}, we find
\begin{th}
If $p>2$, $n>2$, the new stress-energy d-tensors $\tilde{\cal T}_{\alpha\beta}$,
$\tilde{\cal T}_{ij}$ and $\tilde{\cal T}^{(\alpha)(\beta)}_{(i)(j)}$ must
verify the following conservation laws:
\begin{equation}
\left\{\begin{array}{l}\medskip
\displaystyle{\tilde{\cal T}^\mu_{\beta/\mu}+{1\over 2-n}\tilde{\cal  T}_{M/\beta}
+{1\over 2-pn}\tilde{\cal T}_{v/\beta}=-R^m_{\beta\vert m}-P^{(m)}_{(\mu)\beta}
\vert^{(\mu)}_{(m)}}\\\medskip
\displaystyle{{1\over 2-p}\tilde{\cal T}_{T\vert j}+\tilde{\cal  T}^m_{j\vert m}
+{1\over 2-pn}\tilde{\cal T}_{v\vert j}=-P^{(m)}_{(\mu)j}\vert^{(\mu)}_{(m)}}\\
\displaystyle{{1\over 2-p}\tilde{\cal T}_T\vert^{(\alpha)}_{(i)}+{1\over 2-n}
\tilde{\cal  T}_M\vert^{(\alpha)}_{(i)}+\tilde{\cal T}^{(m)(\alpha)}_{(\mu)(i)}
\vert^{(\mu)}_{(m)}=-P^{m(\alpha)}_{\;\;\;\;i\vert m}},
\end{array}\right.
\end{equation}
\end{th}
where
\begin{equation}
\begin{array}{ll}\medskip
R^i_\beta=e^{-2\sigma}\varphi^{im}R_{m\beta},&P^{(i)}_{(\alpha)\beta}=e^{-2\sigma}
\varphi^{im}h_{\alpha\mu}P^{(\mu)}_{(m)\beta},\\
P^{i(\beta)}_{\;(j)}=e^{-2\sigma}\varphi^{im}P_{m(j)}^{\;(\beta)},&
P^i_{(\alpha)j}=e^{-2\sigma}\varphi^{im}h_{\alpha\mu}P^{(\mu)}_{(m)j}.
\end{array}
\end{equation}
\addtocounter{rem}{1}
{\bf Remark \therem} If the conformal function $\sigma$ is independent of partial
directions $x^i_\alpha$, in other words, all functions $\sigma^{(i)}_{(\alpha)}$
vanish, then the conservation laws take the  following simple form:
\begin{equation}
\tilde{\cal T}^\mu_{\beta/\mu}=0,\quad \tilde{\cal T}^m_{i\vert m}=0,\quad
\tilde{\cal T}^{(m)(\alpha)}_{(\mu)(i)}\vert^{(\mu)}_{(m)}=0.
\end{equation}

\section{Maxwell equations of electromagnetic field}

\setcounter{equation}{0}
\hspace{5mm} In order to develope the electromagnetic theory on the generalized
metrical multi-time Lagrange space $GRGML^n_p$, let us consider the {\it canonical Liouville
d-tensor} {\bf C}=$\displaystyle{x^i_\alpha{\partial\over\partial  x^i_\alpha}}$
on $J^1(T,M)$. Using the Cartan canonical connection $C\Gamma$ of $GRGML^n_p$,
we construct the {\it metrical deflection d-tensors} \cite{10}
\begin{equation}
\begin{array}{l}\medskip
\bar D^{(\alpha)}_{(i)\beta}=\left[G^{(\alpha)(\mu)}_{(i)(m)}x^m_\mu\right]_
{/\beta},\\\medskip
D^{(\alpha)}_{(i)j}=\left[G^{(\alpha)(\mu)}_{(i)(m)}x^m_\mu\right]_{\vert j},
\\\medskip
d^{(\alpha)(\beta)}_{(i)(j)}=\left[G^{(\alpha)(\mu)}_{(i)(m)}x^m_\mu\right]
\vert^{(\beta)}_{(j)},
\end{array}
\end{equation}
where $G^{(\alpha)(\beta)}_{(i)(j)}=h^{\alpha\beta}e^{2\sigma}\varphi_{ij}$
is the vertical fundamental  metrical d-tensor  of $GRGML^n_p$ and
$"_{/\beta}"$, $"_{\vert j}"$ or $"\vert^{(\beta)}_{(j)}"$, are the local
covariant derivatives  of $C\Gamma$.
\medskip\\

Taking into account the expressions of the local covariant derivatives of
the Cartan canonical connection $C\Gamma$, we find
\begin{prop}
The metrical  deflection d-tensors of the space $GRGML^n_p$ are given by  the following
formulas,
\begin{equation}
\begin{array}{l}\medskip
\bar D^{(\alpha)}_{(i)\beta}=e^{2\sigma}h^{\alpha\mu}\varphi_{im}\sigma_\beta
x^m_\mu,\\\medskip
D^{(\alpha)}_{(i)j}=-e^{2\sigma}h^{\alpha\mu}\left[\sigma_j\varphi_{im}-
\sigma_i\varphi_{jm}+\sigma_m\varphi_{ij}\right]x^m_\mu,
\\\medskip
d^{(\alpha)(\beta)}_{(i)(j)}=e^{2\sigma}\left[h^{\alpha\beta}\varphi_{ij}+
h^{\alpha\mu}\left(\sigma^{(\beta)}_{(j)}\varphi_{im}-\sigma^{(\beta)}_{(i)}
\varphi_{jm}+\sigma^{(\beta)}_{(m)}\varphi_{ij}\right)x^m_\mu\right].
\end{array}
\end{equation}
\end{prop}
\addtocounter{defin}{1}
{\bf Definition \thedefin} The distinguished 2-form on $J^1(T,M)$,
\begin{equation}
F=F^{(\alpha)}_{(i)j}\delta x^i_\alpha\wedge dx^j+f^{(\alpha)(\beta)}_{(i)(j)}
\delta x^i_\alpha\wedge\delta x^j_\beta,
\end{equation}
where
$F^{(\alpha)}_{(i)j}=\displaystyle{{1\over 2}\left[D^{(\alpha)}
_{(i)j}-D^{(\alpha)}_{(j)i}\right]}$ and
$f^{(\alpha)(\beta)}_{(i)(j)}=\displaystyle{{1\over 2}\left[d^{(\alpha)(\beta)}
_{(i)(j)}-d^{(\alpha)(\beta)}_{(j)(i)}\right]}$, is called the distinguished
{\it electromagnetic 2-form} of the generalized  metrical multi-time Lagrange
space $GRGML^n_p$.
\begin{prop}
The local electromagnetic d-tensors of $GRGML^n_p$ have the expressions,
\begin{equation}
\left\{\begin{array}{l}\medskip
\displaystyle{F^{(\alpha)}_{(i)j}=e^{2\sigma}h^{\alpha\mu}[\varphi_{jm}\sigma_i
-\varphi_{im}\sigma_j]x^m_\mu},\\
\displaystyle{f^{(\alpha)(\beta)}_{(i)(j)}=e^{2\sigma}h^{\alpha\mu}[\varphi_
{im}\sigma_{(j)}^{(\beta)}-\varphi_{jm}\sigma_{(i)}^{(\beta)}]x^m_\mu.}
\end{array}\right.
\end{equation}
\end{prop}

Particularizing the Maxwell equations of electromagnetic field, described in
the general case of a generalized metrical multi-time Lagrange space \cite{10},
we  deduce the main result of the electromagnetism on $GRGML^n_p$.
\begin{th}
The electromagnetic components $F^{(\alpha)}_{(i)j}$ and $f^{(\alpha)(\beta)}_
{(i)(j)}$ of the generalized metrical multi-time Lagrange space $GRGML^n_p$ are
governed by the Maxwell equations:
\medskip
\begin{equation}
\left\{\begin{array}{l}\medskip
\displaystyle{F^{(\alpha)}_{(i)k/\beta}=F_{(i)k}^{(\alpha)}\sigma_\beta+
x^{(\alpha)}_{(i)}\sigma_{\beta\vert k}-x^{(\alpha)}_{(k)}\sigma_{\beta\vert
i}}\\\medskip
\displaystyle{f^{(\alpha)(\gamma)}_{(i)(k)/\beta}=2f^{(\alpha)(\gamma)}_{(i)(k)}
\sigma_\beta+x^{(\alpha)}_{(i)}\sigma^{(\gamma)}_{(k)/\beta}-x^{(\alpha)}_
{(k)}\sigma^{(\gamma)}_{(i)/\beta}}\\\medskip
\displaystyle{\sum_{\{i,j,k\}}F^{(\alpha)}_{(i)j\vert k}=-\sum_{\{i,j,k\}}
x^{(\alpha)}_{(i)}r^s_{mjk}\sigma^{(\mu)}_{(s)}x^m_\mu}\\\medskip
\displaystyle{\sum_{\{i,j,k\}}\left\{F^{(\alpha)}_{(i)j}\vert^{(\gamma)}_{(k)}+
f^{(\alpha)(\gamma)}_{(i)(j)\vert k}\right\}=0}
\\
\displaystyle{\sum_{\{i,j,k\}}f^{(\alpha)(\beta)}_{(i)(j)}\vert^{(\gamma)}_
{(k)}=0,}
\end{array}\right.
\end{equation}
where $x^{(\alpha)}_{(i)}=e^{2\sigma}h^{\alpha\mu}\varphi_{im}x^m_\mu$.
\end{th}
{\bf Proof.} By a direct calculation, the following tensorial identities,
\begin{equation}\label{td}
\begin{array}{l}\medskip
\bar D^{(\alpha)}_{(i)\beta}-x^{(\alpha)}_{(m)}T^m_{\beta i}=2x_{(i)}^{(\alpha)}
\sigma_\beta,\\
d^{(\alpha)(\beta)}_{(i)(j)}+x^{(\alpha)}_{(m)}C^{m(\beta)}_{i(j)}=e^{2\sigma}
h^{\alpha\beta}\varphi_{ij}+2x_{(i)}^{(\alpha)}\sigma_{(j)}^{(\beta)},
\end{array}
\end{equation}
hold  good. The identities  \ref{td}, together with the  general expressions
of Maxwell  equations  for a generalized metrical multi-time  Lagrange space,
imply what we were  looking for. \rule{5pt}{5pt}\medskip\\
\addtocounter{rem}{1}
{\bf  Remark \therem} If the conformal function $\sigma(t^\gamma,x^k,x^k_\gamma)$
is of the form $\sigma(x^k)$, then $f^{(\alpha)(\beta)}_{(i)(j)}$ vanish, and
the Maxwell equations of $GRGML^n_p$ reduce  to
\begin{equation}\label{cme}
\left\{\begin{array}{l}\medskip
\displaystyle{F^{(\alpha)}_{(i)k/\beta}=F_{(i)k}^{(\alpha)}\sigma_\beta+
x^{(\alpha)}_{(i)}\sigma_{\beta\vert k}-x^{(\alpha)}_{(k)}\sigma_{\beta\vert
i}}\\\medskip
\displaystyle{\sum_{\{i,j,k\}}F^{(\alpha)}_{(i)j\vert k}=0}\\
\displaystyle{\sum_{\{i,j,k\}}F^{(\alpha)}_{(i)j}\vert^{(\gamma)}_{(k)}=0}.
\end{array}\right.
\end{equation}

\begin{center}
University POLITEHNICA of Bucharest\\
Department of Mathematics I\\
Splaiul Independentei 313\\
77206 Bucharest, Romania\\
e-mail: mircea@mathem.pub.ro\\
\end{center}

\end{document}